\documentclass{article}    
\usepackage{amstext}
\def\qed{$\rlap{$\sqcap$}\sqcup$}
\usepackage{latexsym}

\begin{document}           

\begin{center}
{\ }\\
{\huge {\bf A non-unimodal codimension 3\\level $h$-vector}} \\ [.250in]
{\large FABRIZIO ZANELLO\\
Dipartimento di Matematica, Universit\`a di Genova, Genova, Italy.\\E-mail: zanello@math.kth.se}
\end{center}

{\ }\\
\\
ABSTRACT. $(1,3,6,10,15,21,28,27,27,28)$ is a level $h$-vector!\\
This example answers negatively the open question as to whether all codimension 3 level $h$-vectors are unimodal.\\
Moreover, using the same (simple) technique, we are able to construct level algebras of codimension 3 whose $h$-vectors have exactly $N$ \lq \lq maxima", for any positive integer $N$.\\
These non-unimodal $h$-vectors, in particular, provide examples of codimension 3 level algebras not enjoying the Weak Lefschetz Property (WLP). Their existence was also an open problem before.\\
In the second part of the paper we further investigate this fundamental property, and show that there even exist codimension 3 level algebras of type 3 without the WLP.\\

{\large

{\ }\\
\\\indent
An important problem in commutative algebra is the study of the properties of level algebras and of their $h$-vectors. This is largely due to the relevance these special algebras have in several other areas of mathematics, such as invariant theory, geometry, combinatorics, etc. (see, e.g., $[GHMS]$ for a broad overview of level algebras).\\\indent
Let us first of all introduce the main definitions we will need in the present paper. We consider standard graded artinian algebras $A=R/I$, where $R=k[x_1,...,x_r]$, $I$ is a homogeneous ideal of $R$, $k$ is a field of characteristic zero and the $x_i$'s all have degree 1.\\\indent
The {\it $h$-vector} of $A$ is $h(A)=h=(h_0,h_1,...,h_e)$, where $h_i=\dim_k A_i$ and $e$ is the last index such that $\dim_k A_e>0$. Since we may suppose, without loss of generality, that $I$ does not contain non-zero forms of degree 1, $r=h_1$ is defined as the {\it codimension} of $A$.\\\indent 
The {\it socle} of $A$ is the annihilator of the maximal homogeneous ideal $\overline{m}=(\overline{x_1},...,\overline{x_r})\subseteq A$, namely $soc(A)=\lbrace a\in A {\ } \mid {\ } a\overline{m}=0\rbrace $. Since $soc(A)$ is a homogeneous ideal, we define the {\it socle-vector} of $A$ as $s(A)=s=(s_0,s_1,...,s_e)$, where $s_i=\dim_k soc(A)_i$. Note that $h_0=1$, $s_0=0$ and $s_e=h_e>0$. The integer $e$ is called the {\it socle degree} of $A$ (or of $h$). The {\it type} of the socle-vector $s$ (or of the algebra $A$) is type$(s)=\sum_{i=0}^es_i$.\\\indent 
If $s=(0,0,...,0,s_e=t)$, we say that the graded algebra $A$ is {\it level} (of type $t$). If, moreover, $t=1$, then $A$ is {\it Gorenstein}. With a slight abuse of notation, we will refer to an $h$-vector as level (or Gorenstein) if it is the $h$-vector of a level (or Gorenstein) algebra.\\\indent
Let us now recall the main facts of the theory of {\it inverse systems}, or {\it Macaulay duality}, which will be a fundamental tool in this note. For a complete introduction, we refer the reader to $[Ge]$ and $[IK]$.\\\indent 
Let $S=k[y_1,...,y_r]$, and consider $S$ as a graded $R$-module where the action of $x_i$ on $S$ is partial differentiation with respect to $y_i$.\\\indent 
There is a one-to-one correspondence between graded artinian algebras $R/I$ and finitely generated graded $R$-submodules $M$ of $S$, where $I=Ann(M)$ is the annihilator of $M$ in $R$ and, conversely, $M=I^{-1}$ is the $R$-submodule of $S$ which is annihilated by $I$ (cf. $[Ge]$, Remark 1), p. 17).\\\indent 
If $R/I$ has socle-vector $s$, then $M$ is minimally generated by $s_i$ elements of degree $i$, for $i=1,2,...,e$, and the $h$-vector of $R/I$ is given by the number of linearly independent partial derivatives obtained in each degree by differentiating the generators of $M$ (cf. $[Ge]$, Remark 2), p. 17).\\\indent 
In particular, level algebras of type $t$ and socle degree $e$ correspond to $R$-submodules of $S$ minimally generated by $t$ elements of degree $e$.\\
\\\indent
It was initially believed that all Gorenstein $h$-vectors of any codimension might be {\it SI-sequences}, i.e., symmetric and with their first half {\it differentiable} (that is, its first difference is again the $h$-vector of an artinian algebra). Instead, this conjecture is false for $r\geq 5$: in particular, Gorenstein $h$-vectors are not even necessarily {\it unimodal} (that is, they no longer increase once they start decreasing) (see $[St2]$, Example 4.3, for the first counter-example, in codimension 13. See also $[BI]$, $[BL]$ and $[Bo1]$).\\\indent
In codimension $r=4$, we do not know whether or not all Gorenstein $h$-vectors are SI-sequences, nor even whether they must be unimodal (see, however, Iarrobino-Srinivasan's paper $[IS]$).\\\indent
Instead, in codimension $r\leq 3$, the conjecture that all Gorenstein $h$-vectors are SI-sequences was proven by Stanley in $[St2]$, Theorem 4.2 (he characterized those $h$-vectors by using a structure theorem of Buchsbaum and Eisenbud, see $[BE]$. See also $[Ma]$ for $r=2$, and our paper $[Za]$, where we show Stanley's result in an elementary fashion).\\\indent
In general, level $h$-vectors were first studied by Stanley (see $[St1]$). Iarrobino, in his 1984 paper $[Ia]$, characterized the level $h$-vectors of codimension 2, while in higher codimensions the situation looks much more difficult and is still mostly unclear. During the last years, however, level $h$-vectors have been extensively investigated; we refer to the memoir $[GHMS]$ for the most comprehensive bibliography on the subject up to 2003.\\
\\\indent
A very natural (and important) problem, especially after Stanley's result on the Gorenstein case, is to determine whether all codimension 3 level $h$-vectors are unimodal (see $[GHMS]$, Question 4.4). In the present note we solve this problem negatively, by using a simple technique which allows us to provide several counter-examples. In particular, we even show that, for any given positive integer $N$, there exist level $h$-vectors of codimension 3 having exactly $N$ \lq \lq maxima".\\
\\\indent
As an easy consequence, we have that not all level algebras of codimension 3 enjoy the {\it Weak Lefschetz Property} ({\it WLP}, briefly). The existence of codimension 3 level algebras without the WLP was also an open problem (see $[GHMS]$, Question 4.4).\\\indent
We recall that an artinian algebra $A=\oplus_{i=0}^e A_i$ is said to enjoy the WLP if there exists a linear form $L$ such that, for all indices $i=0,1,...,e-1$, the multiplication map \lq \lq $\cdot L$" between the $k$-vector spaces $A_i$ and $A_{i+1}$ has maximal rank.\\\indent
The WLP is a fundamental property of artinian algebras and has recently received a lot of attention (e.g., see $[HMNW]$ and $[MM]$ and their bibliographies). In $[HMNW]$, Proposition 3.5, the possible $h$-vectors of algebras with the WLP are characterized: they are unimodal and their increasing part is differentiable. In particular, the above-mentioned Gorenstein algebras of codimension $r\geq 5$ having a non-unimodal $h$-vector do not enjoy the WLP. In codimension 4 (even if, as we said before, no non-unimodal $h$-vector is known), examples of Gorenstein algebras without the WLP have been provided by Ikeda ($[Ik]$, Example 4.4) and Boij ($[Bo2]$, Theorem 3.6).\\\indent
Instead, in codimension 3, so far the only remarkable result in this line of inquiry was that all complete intersections enjoy the WLP ($[HMNW]$, Corollary 2.4). As we mentioned above, until now it was not known if there were level algebras of codimension 3 without the WLP.\\\indent
At the end of this note we further investigate the WLP, supplying several other examples (with a unimodal $h$-vector) of codimension 3 level algebras (some even monomial) without the Weak Lefschetz Property. In particular, we show that there even exist codimension 3 level algebras of type 3 not enjoying the WLP. It remains open if (over a field of characteristic zero) there are Gorenstein or type 2 level algebras of codimension 3 without the WLP. We only mention here that a recent result\\
\\\indent
In order to prove our main result, we need a theorem of Iarrobino, which shows that the partial derivatives of a {\it generic} form (generic in the sense that it corresponds to a point of a suitable non-empty Zariski-open set) \lq \lq intersect" a given inverse system module as little as possible. This result also takes into account the well-known fact that the number of linearly independent $i$-th partial derivatives of a generic form of degree $e$ is the largest possible, i.e. the minimum between ${r-1+e-i \choose e-i}$ and ${r-1+i\choose i}$ (e.g., see $[Ia]$, Proposition 3.4). We should mention that, in the original paper, Iarrobino's result appears with more general hypotheses.\\
\\\indent
{\bf Lemma 1} ($[Ia]$, Theorem 4.8 A). {\it Let $h=(1,h_1,...,h_e)$ be the $h$-vector of a level algebra $A=R/Ann(M)$. Then, if $F$ is a generic form of degree $e$, the level algebra $R/Ann(<M,F>)$ has $h$-vector $H=(1,H_1,...,H_e),$ where, for $i=1,...,e$, $$H_i=\min \lbrace h_i+{r-1+e-i \choose e-i},{r-1+i\choose i}\rbrace .$$}
\\\indent 
{\bf Example 2.} Let us consider a level algebra $A=R/Ann(M)$, where $R=k[x_1,x_2,x_3]$, having $h$-vector $(1,3,6,9,12,15,18,21,24,27)$. For instance, we can construct it by truncating (after degree 9) a Gorenstein algebra having $h$-vector $(1,3,6,9,12,15,18,21,24,27,\\24,21,18,15,12,9,6,3,1)$, which exists by
the above-mentioned characterization of Stanley ($[St2]$, Theorem 4.2). (The fact that the truncation of a Gorenstein algebra, or more generally of a level algebra, is again a level algebra is easy to see and was first noticed by Stanley in $[St1]$.)\\\indent
Then, by Lemma 1, for any generic form $F\in S=k[y_1,y_2,y_3]$ of degree 9, the $h$-vector of the level algebra $R/Ann(<M,F>)$ is $$H=(1,3,6,10,15,21,28,27,27,28).$$
\\\indent
Hence we have shown:\\
\\\indent
{\bf Theorem 3.} {\it There exist non-unimodal level $h$-vectors of codimension 3.}\\
\\\indent
In particular, since the $h$-vectors of algebras with the Weak Lefschetz Property must be unimodal, we immediately have the following:\\
\\\indent
{\bf Corollary 4.} {\it There exist codimension 3 level algebras not enjoying the Weak Lefschetz Property.}\\
\\\indent
{\bf Remark 5.} i). Theorem 3 and Corollary 4 supply a negative answer to all the four sub-questions of $[GHMS]$, Question 4.4.\\\indent
ii). Employing the same idea we used in Example 2, we can construct infinitely many other non-unimodal codimension 3 level $h$-vectors (having larger type and socle degree): for instance, consider a level algebra $A=R/Ann(M)$, having socle degree $e$ large enough, whose $h$-vector ends with an arithmetic progression of the form $(...,t-ip,...,t-2p,t-p,t)$, for any $p\geq 3$ (e.g., $A$ can be constructed as above by truncating a Gorenstein algebra). Then, by Lemma 1, a simple calculation shows that, if we take a generic form $F$ of degree $e$, then the degree $e-p$ entry of the level $h$-vector $H$ of $R/Ann(<M,F>)$ is greater (by 1) than the entry of degree $e-p+1$. Since $t+1>t-p+3$ for $p\geq 3$, it follows that $H$ is not unimodal.\\\indent
For example, applying the above procedure to $(1,3,6,10,15,20,25,30,35,40,45,50,55,\\60,65,70)$ (which is the truncation of the Gorenstein $h$-vector $(1,3,6,...,60,65,70,65,60,...,\\6,3,1)$), we obtain the non-unimodal level $h$-vector $$H=(1,3,6,10,15,21,28,36,45,55,66,65,65,66,68,71).$$\indent
iii). The same technique can be used to show that there exist codimension 3 non-unimodal level $h$-vectors with as many maxima as we want. In particular, there are codimension 3 level $h$-vectors ending with $(...,t,t,t+1,t,t,t+1,...,t,t,t+1)$, where the sequence $t,t,t+1$ may be repeated as many times as we desire.\\\indent
In fact, let us construct a level $h$-vector having exactly $N$ maxima, where $N$ is any positive integer. Consider first a level $h$-vector $h=(1,3,h_2,...,h_e)$, having socle degree $e$ large enough, such that, for some integer $t$, $h_e=t$, $h_{e-1}=t-3$, $h_{e-2}=t-6$, $h_{e-3}=t-9$, $h_{e-4}=t-15$, $h_{e-5}=t-21$, $h_{e-6}=t-27$, $h_{e-7}=t-36$, $h_{e-8}=t-45$, $h_{e-9}=t-54$, and so on, down to $h_{e-3(N-2)-1}=t-{(3(N-2)+1)+2\choose 2}$, $h_{e-3(N-2)-2}=t-{(3(N-2)+1)+2\choose 2}-3(N-1)$, $h_{e-3(N-2)-3}=t-{(3(N-2)+1)+2\choose 2}-6(N-1)$, and $h_{e-3(N-2)-j}={e-3(N-2)-j+2\choose 2}$ for all $j\geq 4$.\\\indent
It is easy to see that such a level $h$-vector actually exists (again, by truncating the suitable Gorenstein $h$-vector of socle degree $2e$, whose existence is clearly guaranteed by $[St2]$, Theorem 4.2). Now, adding a generic form of degree $e$ to an inverse system module generating $h$, by Lemma 1 and a standard computation we obtain the level $h$-vector $$H=(1,3,6,...,H_{e-3(N-2)-4}={(e-3(N-2)-4)+2\choose 2},$$$$H_{e-3(N-2)-3}=t+1,H_{e-3(N-2)-2}=t,H_{e-3(N-2)-1}=t,t+1,t,t,t+1,...,t,t,t+1).$$\indent
For instance, let us construct a codimension 3 level $h$-vector with exactly $N=4$ maxima. Take $e=100$ and $h=(1,3,6,...,h_{90}={90+2\choose 2}=a,a+9,a+18,a+27,a+36,a+42,a+48,a+54,a+57,a+60,h_{100}=a+63=t)$. Using the construction suggested above, we obtain the level $h$-vector
$$H=(1,3,6,...,H_{90}={90+2\choose 2},t+1,t,t,t+1,t,t,t+1,t,t,t+1),$$ which has exactly 4 maxima, as we desired.\\\indent
iv). The existence of non-unimodal level $h$-vectors of codimension 3 immediately implies the existence of non-unimodal level $h$-vectors of codimension $r$, for all $r\geq 4$. In fact, if the level algebra $k[x_1,x_2,x_3]/Ann(M)$ has $h$-vector $h=(1,3,h_2,...,h_e)$, then the level algebra $k[x_1,...,x_r]/Ann(<M,y_4^e,...,y_r^e>)$ has $h$-vector $h^{'}=(1,r,h_2+r-3,....,h_e+r-3)$, which is clearly non-unimodal if $h$ is non-unimodal.\\
\\\indent
Now that we have shown that the property of unimodality in codimension 3 does not extend from Gorenstein to level $h$-vectors of arbitrary type, the main question becomes:\\
\\\indent
{\bf Question 6.} What is the maximum type $t_0$ such that all the codimension 3 level $h$-vectors of type $t\leq t_0$ are unimodal? In particular, is there always unimodality for $t=2$?\\
\\\indent
We only mention here that, as we were finishing to write this paper, Professor Iarrobino has informed us (personal communication) that he and Art Weiss (a graduate student at Tufts University) have just discovered an example of a non-unimodal codimension 3 level $h$-vector of type 5.\\
\\\indent
In the last part of this article we further investigate the Weak Lefschetz Property. We construct other classes of codimension 3 level algebras not enjoying the WLP (these all having a unimodal $h$-vector).\\
\\\indent
{\bf Example 7.} i). Let $M=<y_1^2y_2,y_2^3,y_2^2y_3,y_2y_3^2,y_3^3>\subseteq k[y_1,y_2,y_3]$. A simple computation shows that the $h$-vector of $A=R/I$, where $I=Ann(M)$, is $(1,3,5,5)$.  We want to show that $A$ does not enjoy the Weak Lefschetz Property.\\\indent
Since $I$ is a monomial ideal, it is easy to show that the classes of the monomials of degree $j$ of the inverse system module $M$ (written in the $x_i$'s) generate the vector space $A_j$, for each $j$. Hence, taking the partial derivatives of the generators of $M$, we immediately have that $A_2=<\overline{x_1x_2},\overline{x_1^2},\overline{x_2^2},\overline{x_2x_3},\overline{x_3^2}>$ and
$A_3=<\overline{x_1^2x_2},\overline{x_2^3},\overline{x_2^2x_3},\overline{x_2x_3^2},\overline{x_3^3}>$.\\\indent
Suppose that $A$ enjoys the WLP. Then there exists a linear form $L=ax_1+bx_2+cx_3\in R$ such that the multiplication map \lq \lq $\cdot L$" is a bijection between $A_2$ and $A_3$. We have that $L\cdot \overline{x_1x_2}=\overline{ax_1^2x_2}$ and $L\cdot \overline{x_1^2}=\overline{bx_1^2x_2}$, with $a,b\neq 0$ otherwise \lq \lq $\cdot L$" would not be injective. But the non-zero element $\overline{bx_1x_2-ax_1^2}$ multiplied by $L$ gives $\overline{abx_1^2x_2}-\overline{bax_1^2x_2}=\overline{0}$, a contradiction to the injectivity of \lq \lq $\cdot L$". Hence $A$ does not enjoy the WLP, as we desired.\\\indent
ii). The above example can be generalized to produce codimension 3 monomial level algebras of any socle degree $e\geq 3$ without the WLP. In fact, if we consider the module $$M=<y_1^{e-1}y_2,y_2^e,y_2^{e-1}y_3,y_2^{e-2}y_3^2,...,y_2y_3^{e-1},y_3^e>,$$ then, by the same idea, the algebra $A=R/Ann(M)$, which has $h$-vector $(1,3,5,6,...,e+1,e+2,e+2)$, does not enjoy the WLP, for the map \lq \lq $\cdot L$" between $A_{e-1}$ and $A_e$ is not injective for any linear form $L$.\\
\\\indent
Our main result on level algebras without the WLP is the following:\\
\\\indent
{\bf Proposition 8.} {\it There exist codimension 3 level algebras of type 3 not enjoying the WLP.}\\
\\\indent
{\bf Proof.} The idea is to consider an inverse system module $M$ generated by two suitable forms of degree $e$ in variables $y_1$ and $y_3$ only and by one generic form of the same degree in variables $y_1$ and $y_2$ only. We present here an example with socle degree $e=7$, but examples with higher values of $e$ can be obtained similarly.\\\indent
Let $$M=<y_1^2y_3^5-y_1y_3^6,y_1^3y_3^4-y_1^5y_3^2,F(y_1,y_2)>,$$ where, for example, we take $F(y_1,y_2)=437y_1^7-232y_1^6y_2-423y_1^5y_2^2-567y_1^4y_2^3-769y_1^3y_2^4+831y_1^2y_2^5-916y_1y_2^6-202y_2^7$. Using the computer program $[CoCoA]$, and in particular a CoCoA program which computes the annihilator of a given inverse system module (written by A. Damiano), we see, in this case, that $I=Ann(M)$ has the generator $x_2x_3$ of degree 2 and two generators of degree 5 of the form $c_1x_1^4x_2+c_2x_1^2x_2^3+c_3x_1x_2^4+c_4x_2^5$ and $c_5x_1^3x_2^2+c_6x_1^2x_2^3+c_7x_1x_2^4+c_8x_2^5$, with the $c_i$'s being non-zero constants of the base field $k$ ($I$ has also three generators of degree 6, but we are not interested in them here). Moreover, the $h$-vector of $A=R/I$ is $$h=(1,3,5,7,9,9,6,3).$$\indent
We want to show that $A$ does not enjoy the WLP.\\\indent
Suppose, by contradiction, that $A$ has the WLP. In particular, for some linear form $L=ax_1+bx_2+cx_3$, \lq \lq $\cdot L$" is a bijection between $A_4$ and $A_5$, which both have dimension 9 as $k$-vector spaces. It is easy to see, from the above description of $I$, that two monomial bases for $A_4$ and $A_5$ are $A_4=<\overline{x_1^4},\overline{x_1^3x_2},\overline{x_1^3x_3},\overline{x_1^2x_2^2},\overline{x_1^2x_3^2},\overline{x_1x_2^3},\overline{x_1x_3^3},\overline{x_2^4},\overline{x_3^4}>$ and $A_5=\\<\overline{x_1^5},\overline{x_1^4x_3},\overline{x_1^3x_3^2},\overline{x_1^2x_2^3},\overline{x_1^2x_3^3},\overline{x_1x_2^4},\overline{x_1x_3^4},\overline{x_2^5},\overline{x_3^5}>$. We have that
$$L \cdot \overline{x_1^3x_2}=a\overline{x_1^4x_2}+b\overline{x_1^3x_2^2}=(-{ac_2\over c_1}-{bc_6\over c_5})\overline{x_1^2x_2^3}+(-{ac_3\over c_1}-{bc_7\over c_5})\overline{x_1x_2^4}+(-{ac_4\over c_1}-{bc_8\over c_5})\overline{x_2^5},$$ 
$$L \cdot \overline{x_1^2x_2^2}=(-{ac_6\over c_5}+b)\overline{x_1^2x_2^3}-{ac_7\over c_5}\overline{x_1x_2^4}-{ac_8\over c_5}\overline{x_2^5},$$
$$L \cdot \overline{x_1x_2^3}=a\overline{x_1^2x_2^3}+b\overline{x_1x_2^4},$$
$$L \cdot \overline{x_2^4}=a\overline{x_1x_2^4}+b\overline{x_2^5}.$$\indent
Thus, the image of four linearly independent elements is contained in the span of only three elements, a contradiction to the injectivity of the map \lq \lq $\cdot L$". Therefore $A$ does not enjoy the WLP, as we wanted to show.{\ }{\ }\qed \\
\\\indent
{\bf Remark 9.} i). There are also other methods to construct codimension 3 level algebras without the WLP (or, in some cases, at least very likely without the WLP!). For instance, if we consider an inverse system module $M$ generated by the last $t$ monomials of degree $e$ (according to the lexicographic order) such that the $h$-vector of $R/Ann(M)$ is $(1,3,4,5,...,t-2,t-1,t)$, and we take a form $F=L_1^e+L_2^e$, with the $L_i$'s generic linear forms, then the type $t+1$ level algebra $A=R/Ann(<M,F>)$ has $h$-vector $(1,3,6,7,...,t,t+1,t+1)$ (this fact can be immediately shown by using the original form of $[Ia]$, Theorem 4.8 A, of which we have stated a particular case in Lemma 1 above). Computations with CoCoA seem to indicate that the algebras $A$ constructed in this way do not enjoy the WLP (very reasonably, there seems to be no bijection between $A_{e-1}$ and $A_e$ for any multiplication map \lq \lq $\cdot L$").\\\indent
What we have just said also suggests that one of the \lq \lq strongest results" of $[GHMS]$ on the Weak Lefschetz Property, namely Proposition 5.24 (which proves that, if we add the $e$-th power of one generic linear form to the inverse system module of a given level algebra of socle degree $e$ with the WLP, we obtain another level algebra with the WLP), is very likely stronger than one expected, because it probably cannot be extended to the sum of the powers of 2 or more generic linear forms. A slight modification of the construction employed in Example 2 already supplies a simple proof that such a generalization cannot be done for 6 or more powers: indeed, if the form $F$ of that example were the sum of the powers of $p\geq 6$ generic linear forms, we would end up with a level algebra having $h$-vector $h=(1,3,6,10,15,21,\max\lbrace 18+p,28\rbrace ,27,27,28)$, and this algebra cannot have the WLP (since $h$ is not differentiable for $p\leq 9$ and non-unimodal for $p\geq 10$).\\\indent
ii). At this point, the next natural step will be to determine whether there exist Gorenstein algebras or type 2 level algebras of codimension 3 without the WLP. This problem seems difficult to attack, especially for the Gorenstein case (as we mentioned in the introduction, so far the most important result in this direction is $[HMNW]$, Corollary 2.4, which shows, with a beautiful geometric argument, that all the complete intersections of codimension 3 enjoy the WLP).\\
\\\indent
{\bf Acknowledgements.} We are greatly indebted to Professor A. Iarrobino. In fact, during a recent email exchange between the two of us, he showed us how to construct examples of level $h$-vectors of the form $(1,3,...,t+1,t,t)$ (in particular, he supplied one ending with $(...,15,14,14)$). Also, using a similar technique (inspired by $[Ia]$ and $[CI]$), Professor Iarrobino provided us with examples of non-unimodal codimension 4 level $h$-vectors. (At that moment, we were wondering if those two types of level $h$-vectors might not exist.) It was a refinement of Iarrobino's idea that led us to show the non-unimodality of some codimension 3 level $h$-vectors.\\\indent
We also warmly thank our former Ph.D. supervisor, Prof. Tony Geramita, for his comments and suggestions on the previous versions of this paper.\\\indent
Furthermore, we wish to thank our friend and colleague Alberto Damiano for writing for us a CoCoA program on inverse systems. It has been a very useful instrument for computing some of the examples provided in this note.\\
\\
\\
\\
{\bf \huge References}\\
\\
$[BI]$ {\ } D. Bernstein and A. Iarrobino: {\it A non-unimodal graded Gorenstein Artin algebra in codimension five}, Comm. in Algebra 20 (1992), No. 8, 2323-2336.\\
$[Bo1]$ {\ } M. Boij: {\it Graded Gorenstein Artin algebras whose Hilbert functions have a large number of valleys}, Comm. in Algebra 23 (1995), No. 1, 97-103.\\
$[Bo2]$ {\ } M. Boij: {\it Components of the space parametrizing graded Gorenstein Artinian algebras with a given Hilbert function}, Pacific J. Math. 187 (1999), 1-11.\\ 
$[BL]$ {\ } M. Boij and D. Laksov: {\it Nonunimodality of graded Gorenstein Artin algebras}, Proc. Amer. Math. Soc. 120 (1994), 1083-1092.\\
$[BE]$ {\ } D.A. Buchsbaum and D. Eisenbud: {\it Algebra structures for finite free resolutions, and some structure theorems for ideals of codimension 3}, Amer. J. Math. 99 (1977), 447-485.\\
$[CI]$ {\ } Y.H. Cho and A. Iarrobino: {\it Hilbert Functions and Level Algebras}, J. of Algebra 241 (2001), 745-758.\\
$[CoCoA]$ A. Capani, G. Niesi and L. Robbiano: {\it CoCoA, a system for doing computations in commutative algebra}, available via anonymous ftp, cocoa.dima.unige.it.\\
$[Ge]$ {\ } A.V. Geramita: {\it Inverse Systems of Fat Points: Waring's Problem, Secant Varieties and Veronese Varieties and Parametric Spaces of Gorenstein Ideals}, Queen's Papers in Pure and Applied Mathematics, No. 102, The Curves Seminar at Queen's (1996), Vol. X, 3-114.\\
$[GHMS]$ {\ } A.V. Geramita, T. Harima, J. Migliore and Y.S. Shin: {\it The Hilbert Function of a Level Algebra}, Memoirs of the Amer. Math. Soc., to appear.\\
$[HMNW]$ {\ } T. Harima, J. Migliore, U. Nagel and J. Watanabe: {\it The Weak and Strong Lefschetz Properties for Artinian $K$-Algebras}, J. of Algebra 262 (2003), 99-126.\\
$[Ia]$ {\ } A. Iarrobino: {\it Compressed Algebras: Artin algebras having given socle degrees and maximal length}, Trans. Amer. Math. Soc. 285 (1984), 337-378.\\
$[IK]$ {\ } A. Iarrobino and V. Kanev: {\it Power Sums, Gorenstein Algebras, and Determinantal Loci}, Springer Lecture Notes in Mathematics (1999), No. 1721, Springer, Heidelberg.\\
$[IS]$ {\ } A. Iarrobino and H. Srinivasan: {\it Some Gorenstein Artin algebras of embedding dimension four, I: components of $PGOR(H)$ for $H=(1,4,7,...,1)$}, J. of Pure and Applied Algebra 201 (2005), 62-96.\\
$[Ik]$ {\ } H. Ikeda: {\it Results on Dilworth and Rees numbers of Artinian local rings}, Japan. J. of Math. 22 (1996), 147-158.\\
$[Ma]$ {\ } F.H.S. Macaulay: {\it The Algebraic Theory of Modular Systems}, Cambridge Univ. Press, Cambridge, U.K. (1916).\\
$[MM]$ {\ } J. Migliore and R. Mir\'o-Roig: {\it Ideals of general forms and the ubiquity of the Weak Lefschetz property}, 
J. of Pure and Applied Algebra 182 (2003), 79-107.\\
$[St1]$ {\ } R. Stanley: {\it Cohen-Macaulay Complexes, Higher Combinatorics}, M. Aigner Ed., Reidel, Dordrecht and Boston (1977), 51-62.\\
$[St2]$ {\ } R. Stanley: {\it Hilbert functions of graded algebras}, Adv. Math. 28 (1978), 57-83.\\
$[Za]$ {\ } F. Zanello: {\it Stanley's theorem on codimension 3 Gorenstein $h$-vectors}, Proc. Amer. Math. Soc. 134 (2006), No. 1, 5-8.

}

\end{document}